\documentclass[10pt]{iopart}

\usepackage{iopams}  
\begin{document}

\title[The iteratively regularized Gauss-Newton method]{A convergence analysis of 
the iteratively regularized Gauss-Newton method
under Lipschitz condition}

\newtheorem{theorem}{Theorem}[section]
\newtheorem{lemma}{Lemma}[section]
\newtheorem{corollary}{Corollary}[section]
\newtheorem{proposition}{Proposition}[section]
\def\A{\mathcal A}
\def\B{\mathcal B}
\oddsidemargin 17mm
 \evensidemargin 17mm
 \topmargin 7pt

\catcode`@=11 \@addtoreset{equation}{chapter}
\def\theequation{\thesection.\arabic{equation}}
\catcode`@=12

\author{Qinian Jin}

\address{Department of Mathematics, The University of Texas at Austin, 
Austin, Texas 78712, USA}
\ead{qjin@math.utexas.edu}

\begin{abstract}
In this paper we consider the iteratively regularized Gauss-Newton method 
 for solving nonlinear ill-posed inverse problems. Under merely Lipschitz condition, 
 we prove that this method together with an a posteriori stopping rule defines an 
 order optimal regularization method if the solution is regular in some suitable sense.
\end{abstract}

\maketitle

\section{Introduction}
\setcounter{equation}{0}

In this paper we will consider the nonlinear inverse problems which can be formulated as the 
operator equations
\begin{equation}\label{0}
F(x) = y,
\end{equation}
where $F: D(F)\subset X\to Y$ is a nonlinear operator between the Hilbert spaces $X$ and $Y$ 
with domain $D(F)$ and range $R(F)$. Such problems arise naturally from the parameter 
identification in partial differential equations. For instance, consider the identification 
of the parameter $c$ in the boundary value problem
\begin{equation}\label{01}
\left\{\begin{array}{lll}
-\Delta u+c u=f && \mbox{in } \Omega\\
u=0 && \mbox{on } \partial \Omega
\end{array}\right.
\end{equation}
from the measurement of the state $u$, where $\Omega\subset {R}^n$, $n\le 3$, is a bounded 
domain with smooth boundary $\partial\Omega$, and $f\in L^2(\Omega)$. It is well known that 
(\ref{01}) has a unique solution $u:=u(c)\in H^2(\Omega)\subset L^2(\Omega)$ for each 
\begin{equation*}
c\in D:=\{c\in L^2(\Omega): \|c-\hat{c}\|_{L^2(\Omega)}\le \gamma \mbox{ for some } \hat{c}\ge 0 \mbox{ a.e.}\}
\end{equation*}
with some $\gamma>0$. If we define the operator $F$ as 
\begin{equation*}
F: D\subset L^2(\Omega)\to L^2(\Omega), \quad c\to u(c),
\end{equation*}
then the problem of identifying $c$ is reduced to solving (\ref{0}).

Throghout this paper $\|\cdot\|$ and $(\cdot, \cdot)$ will be used to denote the norms and inner products 
for the both spaces $X$ and $Y$ since there is no confusion. The nonlinear operator $F$ is always assumed to 
be Fr\'echet differentiable, the Fr\'echet derivative of $F$ at $x\in D(F)$ will be denoted as $F'(x)$
and $F'(x)^*$ will be used to denote the adjoint of $F'(x)$. 
We will assume that $y$ is attainable, i.e. $y \in R(F)$. This means that problem 
(\ref{0}) has a solution $x^\dag\in D(F)$, that is  
\begin{equation*}
F(x^\dag)=y.
\end{equation*}
We say problem (\ref{0}) is ill-posed if its solution does not depend continuously on the right
hand side $y$, which is the characteristic property for most of the inverse problems. 
Since the right hand side is usually obtained by measurement, thus, instead of $y$ itself, 
the available data is an approximation $y^\delta$ satisfying
\begin{equation}
\|y^\delta-y\|\le \delta
\end{equation}
with a given small noise level $\delta> 0$. Then the computation of a stable solution of (\ref{0})
from $y^\delta$ becomes an important issue of ill-posed problems, and the regularization
techniques have to be taken into account. 

Tikhonov regularization is one of the well-known method that has been studied extensively in recent years. 
Several a posteriori rules have been suggested to choose the regularization parameter. 
Besides the Morozov's discrepancy principle,  an a posteriori rule has been proposed in
\cite{SEK} to yield optimal rates of convergence. It has been shown in \cite{JH99, TJ} that 
under some reasonable conditions Tikhonov regularization together with the rule in \cite{SEK} 
indeed is order optimal. Moreover, it has also been proved in \cite{TJ} 
that if the solution satisfies suitable source conditions, the same order optimal result 
is still true under merely Lipschitz condition on $F'$.

Iteration methods are also attractive since they are straightforward to implement for the numerical solution of
nonlinear ill-posed problems. In \cite{B}  Bakushinskii proposed the iteratively regularized Gauss-Newton
method in which the iterated solutions $\{x_k^\delta\}$ are defined successively by 
\begin{equation}\label{2}
\fl x_{k+1}^\delta=x_k^\delta-\left(\alpha_k I+F'(x_k^\delta)^*F'(x_k^\delta)\right)^{-1}\left(F'(x_k^\delta)^*(F(x_k^\delta)-y^\delta)+\alpha_k(x_k^\delta-x_0)\right),
\end{equation}
where $x_0^\delta:=x_0$ is an initial guess of $x^\dag$, and $\{\alpha_k\}$ is a given sequence of numbers such that
\begin{equation}\label{3}
\alpha_k>0, \quad 1\le \frac{\alpha_k}{\alpha_{k+1}}\le r \quad \mbox{and}\quad \lim_{k\rightarrow\infty} \alpha_k=0
\end{equation}
for some constant $r>1$. It has been shown in \cite{B, BNS} that if 
\begin{equation}\label{80}
x_0-x^\dag=(F'(x^\dag)^*F'(x^\dag))^{\nu/2} \omega
\end{equation}
for some $\omega\in X$ and some $0<\nu\le 2$, then for the stopping index $N_\delta$ chosen by the 
a priori rule 
\begin{equation*}
\alpha_{N_\delta}\le \left(\frac{\delta}{\|\omega\|}\right)^{2/(1+\nu)}<\alpha_k, \quad 0\le k<N_\delta
\end{equation*}
there holds the order optimal convergence rate
\begin{equation*}
\|x_{N_\delta}^\delta-x^\dag\|\le C\|\omega\|^{1/(1+\nu)}\delta^{\nu/(1+\nu)}
\end{equation*}
with some constant $C$ independent of $\delta$. This rule, however, depends on the knowledge on the smoothness 
of $x_0-x^\dag$, which is difficult to check in practice. Thus a wrong guess of the smoothness will 
lead to a bad choice of $N_\delta$, and consequently to a bad approximation to $x^\dag$. Therefore, a posteriori 
rules, which use only quantities that arise during calculations,  should be considered to 
choose the stopping index of iteration.

The generalized discrepancy principle 
\begin{equation*}
\|F(x_{n_\delta}^\delta)-y^\delta\|\le \tau \delta <\|F(x_k^\delta)-y^\delta\|, \quad 
0\le k<n_\delta
\end{equation*}
has been considered in \cite{BNS} as an a posteriori rule for choosing the stopping index 
$n_\delta$, where $\tau>1$ is a sufficient large number. Under certain conditions, it has been shown that 
\begin{equation*}
\|x_{n_\delta}^\delta-x^\dag\|\le O(\delta^{\nu/(1+\nu)})
\end{equation*}
if $x_0-x^\dag$ satisfies (\ref{80}) with $0<\nu\le 1$. However, with such $n_\delta$, one cannot expect a 
better rate of convergence than $O(\delta^{1/2})$ even if $x_0-x^\dag$ satisfies (\ref{80}) with $\nu>1$.

In order to prevent such saturation, an alternative a posteriori rule has been suggested in \cite{Jin01}
to choose the stopping index $k_\delta$ such that $k_\delta$ is the first integer satisfying 
\begin{equation}\label{4}
\fl \alpha_{k_\delta}\left( F(x_{k_\delta}^\delta)-y^\delta, 
(\alpha_{k_\delta} I+F'(x_{k_\delta}^\delta)F'(x_{k_\delta}^\delta)^*)^{-1}
(F(x_{k_\delta}^\delta)-y^\delta)\right)\le \tau^2 \delta^2,
\end{equation}
where $\tau>1$ is a large number. It has been shown that if $F$ satisfies the condition that there exists a constant $K_0$ such that 
for each pair $x, z\in B_\rho(x^\dag)$ and $w\in X$ there is an element $h(x, z, w)\in X$ such that 
\begin{equation}\label{91}
\fl (F'(x)-F'(z))w=F'(z) h(x, z, w) \quad \mbox{and} \quad \|h(x, z, w)\|\le K_0\|x-z\|\|w\|,
\end{equation}
where $B_\rho(x^\dag)\subset D(F)$ denotes a ball of radius $\rho>0$ around $x^\dag$, then
\begin{equation}\label{82}
\|x_{k_\delta}^\delta-x^\dag\|\le C\|\omega\|^{1/(1+\nu)}\delta^{\nu/(1+\nu)}
\end{equation}
as long as $x_0-x^\dag$ satisfies (\ref{80}) with $0<\nu\le 2$. Moreover, the result in \cite{Jin01} 
implies the convergence rates under more general source conditions,
thus it even applies to exponentially ill-posed problems. 

The structure condition (\ref{91}) on $F$ indeed 
can be verified for many inverse problems. However, there are still some critical cases in which 
this condition is violated. Thus it would be useful to derive some conclusions under
conditions different from or even weaker than (\ref{91}). Some numerical results reported in 
\cite{Jin01} indicates that the order optimal convergence rate holds 
if $x_0-x^\dag$ is smooth enough even if (\ref{91}) is not valid. In this paper we 
will establish some result under merely Lipschitz condition on $F'$.
We will assume that 
\begin{equation}\label{81}
B_\rho(x^\dag)\subset D(F) \quad \mbox{ for some } \rho> 4\|x_0-x^\dag\|
\end{equation}
and that $F'$ satisfies the Lipschitz condition, i.e. there exists a constant $L$ such that
\begin{equation}\label{1}
\|F'(x)-F'(z)\|\le L\|x-z\| \quad \mbox{for all } x, z\in B_\rho(x^\dag).
\end{equation}
Moreover, we will assume that the nonlinear operator $F$ is properly scaled, i.e. 
\begin{equation}\label{94}
\|F'(x)\|\le \alpha_0^{1/2}, \qquad x\in B_\rho(x^\dag).
\end{equation}
This scaling condition can always be fulfilled by multiplying the both sides of (\ref{0}) by
a sufficiently small constant, which then appears as a relaxation parameter in the
iteratively regularized Gauss-Newton method.

The main result of this paper is the following

\begin{theorem}\label{T1}
Assume that (\ref{81}), (\ref{1}) and (\ref{94}) hold and that 
$k_\delta$ is determined by (\ref{4}) with $\tau>2$. If $x_0-x^\dag=F'(x^\dag)^* v$ for some 
$v\in {\mathcal N}(F'(x^\dag)^*)^\perp\subset Y$ and if $L\|v\|$ is sufficiently small, then
\begin{equation}\label{0.8}
\|x_{k_\delta}^\delta-x^\dag\|\le C\inf \left\{\|x_k-x^\dag\|+\frac{\delta}{\sqrt{\alpha_k}}: k=0, 1, \cdots\right\}.
\end{equation}
In particular, if, in addition, $x_0-x^\dag=(F'(x^\dag)^*F'(x^\dag))^{\nu/2}\omega$ for some $\omega\in X$ and some 
$1\le \nu\le 2$, then
\begin{equation}\label{0.8.5}
\|x_{k_\delta}^\delta-x^\dag\|\le C\|\omega\|^{1/(1+\nu)}\delta^{\nu/(1+\nu)},
\end{equation}
where $C$ is a constant independent of $\delta$.
\end{theorem}

We remark that the $\{x_k\}$ appearing in Theorem \ref{T1} is the sequence 
defined by the iteratively regularized Gauss-Newton method (\ref{2}) corresponding to the 
noise-free  case, that is $\{x_k\}$ is defined successively by 
\begin{equation*}
x_{k+1}=x_k-(\alpha_k I+\A_k)^{-1} \left[F'(x_k)^*(F(x_k)-y)+\alpha_k (x_k-x_0)\right].
\end{equation*}
We also remark that we will not specify the smallness of $L\|v\|$ in the above result, 
one can keep track the proof to get an upper bound. The proof will make use of some
ideas developed in \cite{JH99, Jin01, TJ}.

For ease of exposition, throughout this paper we will use the convention 
$A\lesssim B$ to mean that $A\le C B$ for some universal constant $C$ 
independent of $\delta$ and $k$. Thus $A\approx B$ is used to mean 
that $A\lesssim B$ and $B\lesssim A$.

\section{\bf A crucial estimate on convergence rates}
\setcounter{equation}{0}

Throughout this paper, we will use the notations
\begin{eqnarray*}
\fl \quad\quad  \A:=F'(x^\dag)^*F'(x^\dag), & \quad \A_k:=F'(x_k)^*F'(x_k), & \quad  \A_k^\delta:=F'(x_k^\delta)^* F'(x_k^\delta),\\
\fl \quad \quad \B:=F'(x^\dag)F'(x^\dag)^*, & \quad \B_k:=F'(x_k)F'(x_k)^*, & \quad \B_k^\delta:=F'(x_k^\delta) F'(x_k^\delta)^*.
\end{eqnarray*}
We will also use the notations
\begin{eqnarray*}
e_k:=x_k-x^\dag, \qquad  e_k^\delta:=x_k^\delta-x^\dag.
\end{eqnarray*}

The following elementary result will be used frequently.

\begin{lemma}\label{L1}
Let $\{p_k\}$ be a sequence of positive numbers satisfying $\frac{p_k}{p_{k+1}}\le p$ with 
a constant $p>0$. Suppose that the sequence $\{\eta_k\}$ has the property 
$\eta_{k+1}\le p_k+\varepsilon \eta_k$ for all $k$. If $\varepsilon p<1$ and $\eta_0\le \frac{p}{1-\varepsilon p} p_0$, 
then $\eta_k\le \frac{p}{1-\varepsilon p}p_k$ for all $k$.
\end{lemma}

The purpose of this section is to prove the following convergence rate result.

\begin{proposition}\label{P1}
Assume that (\ref{81}), (\ref{1}) and (\ref{94}) hold and that $x_0-x^\dag=F'(x^\dag)^* v$ 
for some $v\in {\mathcal N}(F'(x^\dag)^*)^\perp$. Let $k_\delta$ be the integer determined 
by (\ref{4}) with $\tau>2$. If $L\|v\|$ is sufficiently small, then
\begin{equation*}
\|x_{k_\delta}^\delta-x^\dag\|\le C\|v\|^{1/2}\delta^{1/2},
\end{equation*}
where $C$ is a constant independent of $\delta$.
\end{proposition}

In order to prove Proposition \ref{P1}, we need to give an upper bound on 
$k_\delta$. To this end,  we introduce the integer $\tilde{k}_\delta$ 
which is defined to be the first integer such that
\begin{equation}\label{3.4}
\alpha_{\tilde{k}_\delta}\le \frac{c_0\delta}{\|v\|}<\alpha_k, \quad 0\le k<\tilde{k}_\delta.
\end{equation}
where $c_0$ is a constant such that $0<c_0<\tau-2$.

\begin{lemma}\label{L3}
Under the conditions in Proposition \ref{P1}, if $L\|v\|$ is sufficently small,
then for all integers $0\le k\le \tilde{k}_\delta$  there hold
\begin{equation}\label{3.5}
x_k^\delta\in B_\rho(x^\dag) \qquad \mbox{and}\qquad 
\|e_k^\delta\|\lesssim \alpha_k^{1/2} \|v\|.
\end{equation}
Moreover, for the integer $k_\delta$ determined by (\ref{4}) with $\tau>2$ there holds
$k_\delta\le \tilde{k}_\delta$.
\end{lemma}

\noindent{\bf Proof.}
We first show that $x_k^\delta\in B_\rho(x^\dag)$ for all integers $0\le k\le \tilde{k}_\delta$.
It is clear from (\ref{81}) that this is true for $k=0$. Now for any fixed integer $0<l\le \tilde{k}_\delta$, we assume 
that $x_k^\delta\in B_\rho(x^\dag)$ for all $0\le k<l$ and we are going to show $x_l\in B_\rho(x^\dag)$.
To this end, from the definition of $\{x_k^\delta\}$ it follows that
\begin{eqnarray}\label{6}
\fl e_{k+1}^\delta=\alpha_k(\alpha_k I+\A_k^\delta)^{-1} e_0-(\alpha_k I+\A_k^\delta)^{-1} F'(x_k^\delta)^* 
\left(F(x_k^\delta)-y^\delta-F'(x_k^\delta) e_k^\delta\right).
\end{eqnarray}
Using the condition $e_0=F'(x^\dag)^* v$, we can write
\begin{eqnarray*}
\fl e_{k+1}^\delta=\alpha_k(\alpha_k I+\A)^{-1} e_0 
+\alpha_k\left[(\alpha_k I+\A_k^\delta)^{-1}-(\alpha_k I+\A)^{-1} \right] F'(x^\dag)^* v\\
 -(\alpha_k I+\A_k^\delta)^{-1} F'(x_k^\delta)^*
\left(F(x_k^\delta)-y^\delta-F'(x_k^\delta) e_k^\delta\right).
\end{eqnarray*}
Thus 
\begin{eqnarray}\label{3.1}
\fl \|e_{k+1}^\delta-\alpha_k(\alpha_k I+\A)^{-1} e_0\|
 \le\left\|\alpha_k\left[(\alpha_k I+\A_k^\delta)^{-1}
-(\alpha_k I+\A)^{-1}\right] F'(x^\dag)^* v\right\| \nonumber \\
 +\frac{1}{2\sqrt{\alpha_k}} \|F(x_k^\delta)-y^\delta-F'(x_k^\delta) e_k^\delta\|.
\end{eqnarray}
By using the Lipschitz condition (\ref{1}) we have
\begin{equation}\label{3.2}
\|F(x_k^\delta)-y^\delta-F'(x_k^\delta) e_k^\delta\|\le \delta+ \frac{1}{2} L\|e_k^\delta\|^2.
\end{equation}
Moreover, note that
\begin{eqnarray*}
\fl \alpha_k\left[(\alpha_k I+\A_k^\delta)^{-1}-(\alpha_k I+\A)^{-1}\right] F'(x^\dag)^* v\\
 =\alpha_k (\alpha_k I +\A_k^\delta)^{-1} (\A-\A_k^\delta)(\alpha_k I+\A)^{-1} F'(x^\dag)^* v\\
 =\alpha_k (\alpha_k I +\A_k^\delta)^{-1} F'(x_k^\delta)^*(F'(x^\dag)
-F'(x_k^\delta)) (\alpha_k I+\A)^{-1} F'(x^\dag)^* v\\
 +\alpha_k(\alpha_k I+\A_k^\delta)^{-1} \left(F'(x^\dag)^*-F'(x_k^\delta)^*\right) 
F'(x^\dag)(\alpha_k I+\A)^{-1} F'(x^\dag)^* v
\end{eqnarray*}
We then use (\ref{1}) to obtain
\begin{equation}\label{3.3}
\left\|\alpha_k\left[(\alpha_k I+\A_k^\delta)^{-1}
-(\alpha_k I+\A)^{-1}\right] F'(x^\dag)^* v\right\|
\le 2 L\|v\| \|e_k^\delta\|.
\end{equation}
Combining (\ref{3.1}), (\ref{3.2}) and (\ref{3.3}) gives
\begin{equation*}
\|e_{k+1}^\delta-\alpha_k (\alpha_k I+\A)^{-1} e_0\|\le 2 L\|v\| \|e_k^\delta\|
+\frac{1}{4\sqrt{\alpha_k}}L \|e_k^\delta\|^2+\frac{\delta}{2\sqrt{\alpha_k}}.
\end{equation*}
Let
\begin{equation}\label{3.61}
\beta_k:=\|\alpha_k(\alpha_k I+\A)^{-1} e_0\|.
\end{equation}
Then we have
\begin{equation}\label{3.62}
\|e_{k+1}^\delta\|\le \beta_k+2 L\|v\|\|e_k^\delta\| 
+\frac{1}{4\sqrt{\alpha_k}} L\|e_k^\delta\|^2 +\frac{\delta}{2\sqrt{\alpha_k}}.
\end{equation}
Note that for $0\le k<\tilde{k}_\delta$ we have 
$\frac{\delta}{2\sqrt{\alpha_k}}\le \frac{1}{2c_0} \alpha_k^{1/2} \|v\|$; 
note also that $\beta_k\le \frac{1}{2}\alpha_k^{1/2}\|v\|$. We thus obtain
\begin{equation}\label{3.63}
\|e_{k+1}^\delta\|\le \left(\frac{1}{2}+\frac{1}{2 c_0}\right) \alpha_k^{1/2}\|v\|
+2 L\|v\|\|e_k^\delta\| +\frac{1}{4\sqrt{\alpha_k}} L\|e_k^\delta\|^2. 
\end{equation}
This and (\ref{3}) imply
\begin{equation*}
\frac{\|e_{k+1}^\delta\|}{\sqrt{\alpha_{k+1}}}\le r^{1/2}
\left[\left(\frac{1}{2}+\frac{1}{2c_0}\right)\|v\| 
+2 L\|v\|\frac{\|e_k^\delta\|}{\sqrt{\alpha_k}} 
+\frac{1}{4} L\left(\frac{\|e_k^\delta\|}{\sqrt{\alpha_k}}\right)^2\right].
\end{equation*}
By induction, (\ref{94}) and $e_0=F'(x^\dag)^* v$ we can show that 
if $L\|v\|$ is sufficiently small then
\begin{equation}\label{3.64}
\frac{\|e_k^\delta\|}{\sqrt{\alpha_k}}\le  r^{1/2} \left(1+\frac{1}{c_0}\right) \|v\| \quad \mbox{ for } 0\le k\le l.
\end{equation}
Combining this with (\ref{3.62}) and noting that 
$\frac{\delta}{2\sqrt{\alpha_k}}\le \frac{1}{2\sqrt{c_0}} \|v\|^{1/2}\delta^{1/2}$, 
we have for $0\le k<l$ that
\begin{equation*}
\|e_{k+1}^\delta\|\le \beta_k+\frac{1}{2\sqrt{c_0}}\|v\|^{1/2}\delta^{1/2}+CL\|v\|\|e_k^\delta\|,
\end{equation*}
where here and below $C$ denotes a universal constant independent of $\delta$ and $k$. 
Recall that $\beta_k\le r \beta_{k+1}$ which was proved in \cite[Lemma 3.4]{Jin01}.
By the smallness of $L\|v\|$, we may apply Lemma \ref{L1} to conclude
\begin{equation*}
\|e_k^\delta\|\le 2\beta_k+\frac{1}{\sqrt{c_0}}\|v\|^{1/2} \delta^{1/2},
\end{equation*}
since, due to (\ref{94}), this is true for $k=0$. Note that 
$\beta_k\le \|e_0\|\le \frac{\rho}{4}$, the above inequality implies 
that $x_k^\delta\in B_\rho(x^\dag)$ for all $0\le k\le l$. 
We thus obtain $x_k^\delta\in B_\rho(x^\dag)$ for all $0\le k\le \tilde{k}_\delta$. 
In the meanwhile, (\ref{3.64}) gives the desired estimates
in (\ref{3.5}). 

In order to prove $k_\delta\le \tilde{k}_\delta$, we note that the combination of 
(\ref{3.63}) and (\ref{3.64}) gives
\begin{equation*}
\|e_{k+1}^\delta\|\le \left(\frac{1}{2}+\frac{1}{2 c_0}\right) \alpha_k^{1/2}\|v\|+CL\|v\| \|e_k^\delta\|.
\end{equation*}
Thus, by Lemma \ref{L1},
\begin{equation*}
\|e_k^\delta\|\le \left(1+\frac{1}{c_0}\right)\alpha_k^{1/2} \|v\| \quad \mbox{ for all } 0\le k\le \tilde{k}_\delta.
\end{equation*}
This together with (\ref{1}) and (\ref{3.4}) implies that
\begin{eqnarray*}
\fl \left\|\alpha_{\tilde{k}_\delta}^{1/2} (\alpha_{\tilde{k}_\delta} I+\B_{\tilde{k}_\delta}^\delta)^{-1/2}
\left(F(x_{\tilde{k}_\delta}^\delta)-y^\delta\right)\right\|\\
\fl \qquad\qquad \quad  \le \delta+ \left\|\alpha_{\tilde{k}_\delta}^{1/2} (\alpha_{\tilde{k}_\delta} I+\B_{\tilde{k}_\delta}^\delta)^{-1/2}
F'(x_{\tilde{k}_\delta}^\delta) e_{\tilde{k}_\delta}^\delta\right\|
+\left\|F(x_{\tilde{k}_\delta}^\delta)-y-F'(x_{\tilde{k}_\delta}^\delta) e_{\tilde{k}_\delta}^\delta\right\|\\
 \fl \qquad \qquad \quad \le \delta +\alpha_{\tilde{k}_\delta}^{1/2}\|e_{\tilde{k}_\delta}^\delta\|
+\frac{1}{2} L\|e_{\tilde{k}_\delta}^\delta\|^2\\
\fl \qquad\qquad\quad  \le \delta+ \left(1+\frac{1}{c_0}\right) \|v\| \alpha_{\tilde{k}_\delta} 
+\frac{1}{2} \left(1+\frac{1}{c_0}\right)^2 L\|v\|^2\alpha_{\tilde{k}_\delta}\\
\fl \qquad\qquad\quad  \le \delta +\left(1+\frac{1}{c_0}\right)c_0 \delta 
+\frac{1}{2} \left(1+\frac{1}{c_0}\right)^2 c_0 L\|v\| \delta.
\end{eqnarray*}
Recall that $\tau>2$ and $0<c_0<\tau-2$, we can see for sufficiently small 
$L\|v\|$ there holds
\begin{equation*}
\left\|\alpha_{\tilde{k}_\delta}^{1/2} (\alpha_{\tilde{k}_\delta} I+\B_{\tilde{k}_\delta}^\delta)^{-1/2}
\left(F(x_{\tilde{k}_\delta}^\delta)-y^\delta\right)\right\|\le \tau \delta.
\end{equation*}
By the definition of $k_\delta$, we thus conclude that $k_\delta\le \tilde{k}_\delta$. \hfill  $\Box$\\

Now we are ready to give the proof of Proposition \ref{P1}.\\

\noindent{\bf Proof of Proposition \ref{P1}.} 
If $k_\delta=0$, then, by the definition of $k_\delta$ and (\ref{94}), we have
\begin{equation*}
\|F(x_0)-y^\delta\|\lesssim \delta.
\end{equation*}
This together with (\ref{1}) gives
\begin{eqnarray*}
\|F'(x^\dag) e_0\|&\le \|F(x_0)-y-F'(x^\dag) e_0\|+\|F(x_0)-y\|\\
&\lesssim \frac{1}{2}L\|e_0\|^2 +\delta.
\end{eqnarray*}
Thus, by using $e_0=F'(x^\dag)^* v$, it follows 
\begin{eqnarray*}
\|e_0\|&=(e_0, F'(x^\dag)^*v)^{1/2}=(F'(x^\dag) e_0, v)^{1/2}
\le \|F'(x^\dag) e_0\|^{1/2}\|v\|^{1/2}\\
&\lesssim \sqrt{L\|v\|} \|e_0\|+\|v\|^{1/2}\delta^{1/2}.
\end{eqnarray*}
By the smallness of $L\|v\|$, we obtain
\begin{equation*}
\|e_{k_\delta}^\delta\|=\|e_0\|\lesssim \|v\|^{1/2}\delta^{1/2}.
\end{equation*}

Therefore we may assume $k_\delta>0$ in the following.  Recall that in the proof 
of Lemma \ref{L3} we have obtained the following two estimates
\begin{equation}\label{3.5.5}
\|e_{k+1}^\delta\|\lesssim \beta_k+\|v\|^{1/2}\delta^{1/2}+L\|v\| \|e_k^\delta\|
\end{equation}
and 
\begin{equation}\label{3.6}
\|e_k^\delta\|\lesssim \beta_k+\|v\|^{1/2}\delta^{1/2} 
\end{equation}
for all $0\le k< \tilde{k}_\delta$. 

Now we set  
\begin{equation*}
\beta_k^\delta:=\|\alpha_k(\alpha_k I+\A_k^\delta)^{-1} e_0\|.
\end{equation*}
Then it follows from (\ref{3.3}) that  
\begin{equation*}
|\beta_k-\beta_k^\delta|\le \|\alpha_k\left[(\alpha_k I+\A_k^\delta)^{-1} -(\alpha_k I+\A)^{-1}\right] e_0\|
\le 2 L\|v\|\|e_k^\delta\|.
\end{equation*}
This together with (\ref{3.5.5}) and (\ref{3.6}) implies for small $L\|v\|$ that
\begin{eqnarray}
\|e_k^\delta\|\lesssim \beta_k^\delta+\|v\|^{1/2} \delta^{1/2},\label{3.8.1}\\
\|e_{k+1}^\delta\|\lesssim \beta_k^\delta+\|v\|^{1/2} \delta^{1/2}.\label{3.8.2}
\end{eqnarray}
We need to estimate $\beta_k^\delta$. We first have
\begin{eqnarray*}
\fl \left(\beta_k^\delta\right)^2 =\left(\alpha_k (\alpha_k I+\A_k^\delta)^{-1} e_0, 
\alpha_k (\alpha_k I+\A_k^\delta)^{-1} F'(x^\dag)^* v\right)\\
\fl\quad\quad\quad \quad =\left(\alpha_k (\alpha_k I+\A_k^\delta)^{-1} e_0, \alpha_k (\alpha_k I+\A_k^\delta)^{-1}\left[ F'(x_k^\delta)^*+(F'(x^\dag)^*-F'(x_k^\delta)^*)\right]v\right)\\
=\left(\alpha_k^{3/2}(\alpha_k I+\B_k^\delta)^{-3/2} F'(x_k^\delta) e_0, \alpha_k^{1/2} (\alpha_k I+\B_k^\delta)^{-1/2} v\right)\\
 +\left(\alpha_k (\alpha_k I +\A_k^\delta)^{-1} e_0, \alpha_k (\alpha_k I +\A_k^\delta)^{-1} 
\left[F'(x^\dag)^*-F'(x_k^\delta)^*\right] v\right)\\
\le  \gamma_k^\delta \|v\|+\beta_k^\delta L\|v\| \|e_k^\delta\|,
\end{eqnarray*}
where 
\begin{equation*}
\gamma_k^\delta:=\left\|\alpha_k^{3/2}(\alpha_k I+\B_k^\delta)^{-3/2} F'(x_k^\delta) e_0\right\|.
\end{equation*}
Therefore  
\begin{equation}\label{3.9}
\beta_k^\delta\le \sqrt{\gamma_k^\delta} \|v\|^{1/2}+L\|v\|\|e_k^\delta\|.
\end{equation}
In order to estimate $\gamma_k^\delta$, we observe that (\ref{6}) implies 
\begin{eqnarray*}
\fl \alpha_k^{3/2} (\alpha_k I+\B_k^\delta)^{-3/2} F'(x_k^\delta) e_0
=\alpha_k^{1/2} (\alpha_k I+\B_k^\delta)^{-3/2} \B_k^\delta 
\left(F(x_k^\delta)-y^\delta-F'(x_k^\delta) e_k^\delta\right)\\
+\alpha_k^{1/2} (\alpha_k I+\B_k^\delta)^{-1/2} F'(x_k^\delta) e_{k+1}^\delta.
\end{eqnarray*}
Thus
\begin{eqnarray*}
\fl \gamma_k^\delta\le \left\|\alpha_k^{1/2} (\alpha_k I+\B_k^\delta)^{-1/2} F'(x_k^\delta) e_{k+1}^\delta \right\|
+\left\|F(x_k^\delta)-y^\delta-F'(x_k^\delta) e_k^\delta\right\|\\
\fl \qquad\qquad \le \left\|\alpha_k^{1/2} (\alpha_k I+\B_k^\delta)^{-1/2} \left(F(x_{k+1}^\delta)-y^\delta\right)\right\|
+\left\|\left(F'(x_{k+1}^\delta)-F'(x_k^\delta)\right) e_{k+1}^\delta\right\|\\
+\left\|F(x_{k+1}^\delta)-y-F'(x_{k+1}^\delta) e_{k+1}^\delta\right\|
+\left\|F(x_k^\delta)-y-F'(x_k^\delta) e_k^\delta\right\|+2\delta.
\end{eqnarray*}
It then follows from (\ref{1}) that 
\begin{eqnarray*}
\fl \gamma_k^\delta
\le \left\|\alpha_k^{1/2} (\alpha_k I+\B_k^\delta)^{-1/2} \left(F(x_{k+1}^\delta)-y^\delta\right) \right\|
+2\delta +L\|e_k^\delta\|^2 +2 L\|e_{k+1}^\delta\|^2.
\end{eqnarray*}
Using \cite[Proposition 3.4]{TJ}, and noting that
\begin{equation*}
L\|x_{k+1}^\delta-x_k^\delta\|\le L\left(\|e_k^\delta\|+\|e_{k+1}^\delta\|\right)\lesssim L\|v\|\alpha_k^{1/2},
\end{equation*}
we can conclude for sufficiently small $L\|v\|$ there holds
\begin{eqnarray*}
\fl  \left\|\alpha_k^{1/2} (\alpha_k I+\B_k^\delta)^{-1/2} \left(F(x_{k+1}^\delta)-y^\delta\right) \right\|
 \lesssim \left\|\alpha_k^{1/2} 
(\alpha_k I+\B_{k+1}^\delta)^{-1/2}\left( F(x_{k+1}^\delta)-y^\delta\right) \right\|
\end{eqnarray*}
Let $\{E_\lambda\}$ be the spectral family generated by the self-adjoint operator $\B_{k+1}^\delta$. Then, 
by using (\ref{3}), we have 
\begin{eqnarray*}
\fl \left\|\alpha_k^{1/2} 
(\alpha_k I+\B_{k+1}^\delta)^{-1/2}\left( F(x_{k+1}^\delta)-y^\delta\right) \right\|^2
=\int_0^\infty \frac{\alpha_k}{\alpha_k+\lambda} d \|E_\lambda (F(x_{k+1}^\delta)-y^\delta)\|^2\\
\le r\int_0^\infty \frac{\alpha_{k+1}}{\alpha_{k+1}+\lambda} d \|E_\lambda (F(x_{k+1}^\delta)-y^\delta)\|^2\\
=r\left\|\alpha_{k+1}^{1/2} 
(\alpha_{k+1} I+\B_{k+1}^\delta)^{-1/2}\left( F(x_{k+1}^\delta)-y^\delta\right) \right\|^2.
\end{eqnarray*}
Thus 
\begin{eqnarray*}
\fl  \left\|\alpha_k^{1/2} (\alpha_k I+\B_k^\delta)^{-1/2} \left(F(x_{k+1}^\delta)-y^\delta\right) \right\|
  \lesssim \left\|\alpha_{k+1}^{1/2} 
(\alpha_{k+1} I+\B_{k+1}^\delta)^{-1/2}\left( F(x_{k+1}^\delta)-y^\delta\right) \right\|.
\end{eqnarray*}
Therefore 
\begin{equation*}
\fl \gamma_k^\delta\lesssim \left\|\alpha_{k+1}^{1/2} (\alpha_{k+1} I+\B_{k+1}^\delta)^{-1/2}
\left( F(x_{k+1}^\delta)-y^\delta\right) \right\|+\delta +L\|e_k^\delta\|^2 +L\|e_{k+1}^\delta\|^2.
\end{equation*}
This together with (\ref{3.9}) gives
\begin{eqnarray*}
 \beta_k^\delta\lesssim &\left\|\alpha_{k+1}^{1/2} (\alpha_{k+1} I+\B_{k+1}^\delta)^{-1/2}
\left( F(x_{k+1}^\delta)-y^\delta\right) \right\|^{1/2} \|v\|^{1/2} \\
 & +\|v\|^{1/2}\delta^{1/2}+ \sqrt{L\|v\|}\left(\|e_k^\delta\|+\|e_{k+1}^\delta\|\right).
\end{eqnarray*}
Combining this with (\ref{3.8.1}) and (\ref{3.8.2}) yields 
\begin{eqnarray*}
\fl \beta_k^\delta \lesssim \left\|\alpha_{k+1}^{1/2} (\alpha_{k+1} I+\B_{k+1}^\delta)^{-1/2}
\left( F(x_{k+1}^\delta)-y^\delta\right) \right\|^{1/2}\|v\|^{1/2}
+\|v\|^{1/2}\delta^{1/2} +\sqrt{L\|v\|} \beta_k^\delta.
\end{eqnarray*}
Using the smallness of $L\|v\|$ we obtain
\begin{equation*}
\fl \beta_k^\delta \lesssim \left\|\alpha_{k+1}^{1/2} (\alpha_{k+1} I+\B_{k+1}^\delta)^{-1/2}
\left( F(x_{k+1}^\delta)-y^\delta\right) \right\|^{1/2}\|v\|^{1/2}
+\|v\|^{1/2}\delta^{1/2}.
\end{equation*}
It then follows from (\ref{3.8.2}) that for all $0<k\le \tilde{k}_\delta$ there holds
\begin{equation*}
\|e_k^\delta\|\lesssim \left\|\alpha_{k}^{1/2} (\alpha_{k} I+\B_{k}^\delta)^{-1/2}
\left( F(x_{k}^\delta)-y^\delta\right) \right\|^{1/2}\|v\|^{1/2}+\|v\|^{1/2}\delta^{1/2}.
\end{equation*}
Thus, by setting $k=k_\delta$ in the above inequality and using the 
definition of $k_\delta$,  we obtain the desired estimate. \hfill $\Box$

\section{\bf A key inequality }
\setcounter{equation}{0}

The main result of this section is the following inequality.

\begin{proposition}\label{P2}
Assume that (\ref{81}), (\ref{1}) and (\ref{94}) hold and that $x_0-x^\dag=F'(x^\dag)^* v$ for some $v\in 
{\mathcal N}(F'(x^\dag)^*)^\perp$. If $L\|v\|$ is sufficiently small, then for
any integer $k_\delta\le  k\le \tilde{k}_\delta$ there holds
\begin{equation*}
\|e_{k_\delta}\|\lesssim \|e_k\|+\frac{\|\alpha_{k_\delta}^{1/2}
(\alpha_{k_\delta} I+\B)^{-1/2}(F(x_{k_\delta})-y)\|}{\sqrt{\alpha_k}}.
\end{equation*}
\end{proposition}

The proof of this result will employ Proposition \ref{P1} and 
the following two auxiliary results which are of independent interest. 

\begin{lemma}\label{L2}
Under the conditions in Proposition \ref{P2}, if $L\|v\|$ is sufficently small,
then for all $k\ge 0$ there hold
\begin{equation}\label{3.0}
x_k\in B_\rho(x^\dag) \qquad \mbox{and} \qquad \|e_k\|\lesssim \alpha_k^{1/2}\|v\|.
\end{equation}
Moreover, for all integers $0\le k\le l$ there hold
\begin{equation}\label{3.100}
\|e_k\|\approx \beta_k, 
 \quad \|e_{k-1}\|\lesssim \|e_k\|\quad \mbox{and} \quad 
\|e_l\|\lesssim \|e_k\|,
\end{equation}
where $\beta_k$ is defined as in (\ref{3.61}).
\end{lemma}

\noindent{\bf Proof.}
From the definition of $\{x_k\}$ it follows easily that
\begin{equation}\label{3.0.6}
\fl e_{k+1}=\alpha_k(\alpha_k I+\A_k)^{-1} e_0-(\alpha_k I+\A_k)^{-1} F'(x_k)^* \left(F(x_k)-y-F'(x_k) e_k\right).
\end{equation}
Then we can use the similar argument in the proof of Lemma \ref{L3} to conclude (\ref{3.0}) and the estimate
\begin{equation}\label{70}
\left|\|e_{k+1}\|-\beta_k\right|\le \|e_{k+1}-\alpha_k(\alpha_k I+\A)^{-1} e_0\|\lesssim L\|v\|\|e_k\|
\end{equation}
Thus, by Lemma \ref{L1}, we have 
\begin{equation*}
\|e_k\|\le 2\beta_k.
\end{equation*}
Note that $\beta_k$ is non-increasing, we can use 
(\ref{70}) again to obtain 
\begin{equation*}
\|e_{k+1}\|\ge \beta_k-CL\|v\|\|e_k\|\ge (1-CL\|v\|)\beta_k\ge (1-CL\|v\|)\beta_{k+1}.
\end{equation*}
Therefore $\|e_k\|\approx \beta_k$. As an immediate consequence, we have 
\begin{equation*}
\|e_l\|\lesssim \beta_l\le \beta_k\lesssim \|e_k\|
\end{equation*}
for all $0\le k\le l$. In order to show 
\begin{equation*}
\|e_{k-1}\|\lesssim \|e_k\|,
\end{equation*}
it suffices to show 
\begin{equation*}
\beta_{k-1}\lesssim \beta_k. 
\end{equation*}
However, this last inequality has been verified in \cite[Lemma 4.3]{Jin01}.
\hfill $\Box$

\begin{lemma}\label{L5}
Under the conditions in Proposition \ref{P2}, if $L\|v\|$ is sufficiently small, then for
all integers $0\le k\le \tilde{k}_\delta$ there hold
\begin{equation*}
\|x_k^\delta-x_k\|\le \frac{\delta}{\sqrt{\alpha_k}}.
\end{equation*}
\end{lemma}

\noindent{\bf Proof.}
By setting 
\begin{equation}
u_k:=F(x_k)-y-F'(x_k) e_k, \quad u_k^\delta:=F(x_k^\delta)-y -F'(x_k^\delta) e_k^\delta,
\end{equation}
it then follows from (\ref{6}) and (\ref{3.0.6}) that 
\begin{equation*}
x_{k+1}^\delta-x_{k+1}=I_1+I_2+I_3+I_4,
\end{equation*}
where 
\begin{eqnarray*}
I_1&:=\alpha_k \left[(\alpha_k I+\A_k^\delta)^{-1}-(\alpha_k I+\A_k)^{-1}\right] e_0,\\
I_2&:=(\alpha_k I+\A_k^\delta)^{-1} F'(x_k^\delta)^* (y^\delta-y),\\
I_3&:=\left[(\alpha_k I+\A_k)^{-1} F'(x_k)^*-(\alpha_k I+\A_k^\delta)^{-1} F'(x_k^\delta)^*\right] u_k^\delta,\\
I_4&:=(\alpha_k I+\A_k)^{-1} F'(x_k)^* (u_k-u_k^\delta).
\end{eqnarray*}
It is clear that $\|I_2\|\le \frac{\delta}{2\sqrt{\alpha_k}}$. In order to estimate $I_1$, recall that 
$e_0=F'(x^\dag)^* v$, we can write
\begin{eqnarray*}
I_1&=\alpha_k\left[(\alpha_k I+\A_k^\delta)^{-1}-(\alpha_k I+\A_k)^{-1}\right] F'(x_k)^* v\\
&\quad \,+\alpha_k \left[(\alpha_k I+\A_k^\delta)^{-1} -(\alpha_k I+\A_k)^{-1} \right] (F'(x^\dag)^*-F'(x_k)^*) v.
\end{eqnarray*}
We may use (\ref{1}), the similar argument in deriving (\ref{3.3}), and Lemma \ref{L2} to obtain
$$
\|I_1\|\lesssim L\|v\|\|x_k^\delta-x_k\|+\frac{1}{\sqrt{\alpha_k}} L^2 \|v\|\|e_k\|\|x_k^\delta-x_k\|
\lesssim L\|v\| \|x_k^\delta-x_k\|.
$$
Similarly, for $I_3$ we have
\begin{eqnarray*}
\|I_3\|&\lesssim \left\|(\alpha_k I+\A_k^\delta)^{-1} \left[F'(x_k^\delta)^*-F'(x_k)^*\right] u_k^\delta\right\|\\
&\quad\,+ \left\|\left[(\alpha_k I+\A_k^\delta)^{-1}-(\alpha_k I+\A_k)^{-1}\right] F'(x_k)^* u_k^\delta\right\|\\
&\lesssim \frac{1}{\alpha_k} L\|u_k^\delta\|\|x_k^\delta-x_k\|\\
&\lesssim \frac{1}{\alpha_k} L^2\|e_k^\delta\|^2 \|x_k^\delta-x_k\|\\
&\lesssim L\|v\|\|x_k^\delta-x_k\|.
\end{eqnarray*}
By using (\ref{1}), Lemma \ref{L2} and Lemma \ref{L3} we have
\begin{eqnarray*}
\|I_4\|&\le \frac{1}{2\sqrt{\alpha_k}}\|u_k^\delta-u_k\|\\
&\lesssim \frac{1}{\sqrt{\alpha_k}} \left(L\|x_k^\delta-x_k\|^2+L\|e_k^\delta\|\|x_k^\delta-x_k\|\right)\\
&\lesssim \frac{1}{\sqrt{\alpha_k}}L\left(\|e_k^\delta\|+\|e_k\|\right) \|x_k^\delta-x_k\|\\
&\lesssim L\|v\|\|x_k^\delta-x_k\|. 
\end{eqnarray*}
Combining the above estimates on $I_1$, $I_2$, $I_3$ and $I_4$ we conclude that there is a constant $C$
independent of $\delta$ and $k$ such that for all $0\le k\le \tilde{k}_\delta$
\begin{equation*}
\|x_{k+1}^\delta-x_{k+1}\|\le \frac{\delta}{2\sqrt{\alpha_k}}+ C L\|v\|\|x_k^\delta-x_k\|.
\end{equation*}
Since $L\|v\|$ is small, an application of Lemma \ref{L1} gives the desired estimates.
\hfill$\Box$\\

Now we are in a position to prove Proposition \ref{P2}.\\

\noindent{\bf Proof of Proposition \ref{P2}.}
Let 
\begin{equation*}
J:=\left[(\alpha_{k_\delta-1}(\alpha_{k_\delta-1} I+\A)^{-1}-\alpha_{k-1}(\alpha_{k-1}I+\A)^{-1}\right] e_0.
\end{equation*}
Then it follows from (\ref{3.0.6}), (\ref{1}), and Lemma \ref{L2} that 
\begin{eqnarray}\label{3.10}
\fl \|x_{k_\delta}-x_k\|\le \|J\|+\left\|\alpha_{k_\delta-1}
\left[\alpha_{k_\delta-1} I+\A_{k_\delta-1})^{-1}-(\alpha_{k_\delta-1} I+\A)^{-1}\right]
F'(x^\dag)^* v \right\|\nonumber\\
 +\|\alpha_{k-1}\left[(\alpha_{k-1} I+\A_{k-1})^{-1}-(\alpha_{k-1} I+\A)^{-1}\right] F'(x^\dag)^* v\| \nonumber\\
 +\frac{1}{\sqrt{\alpha_{k_\delta-1}}}\|F(x_{k_\delta-1})-y-F'(x_{k_\delta-1})e_{k_\delta-1}\|\nonumber\\
 +\frac{1}{\sqrt{\alpha_{k-1}}}\|F(x_{k-1})-y-F'(x_{k-1})e_{k-1}\|\nonumber\\
\lesssim \|J\|+L\|v\|\left(\|e_{k_\delta-1}\|+\|e_{k-1}\|\right)+\frac{L\|e_{k_\delta-1}\|^2}{\sqrt{\alpha_{k_\delta-1}}}
+\frac{L\|e_{k-1}\|^2}{\sqrt{\alpha_{k-1}}}\nonumber\\
\lesssim \|J\|+L\|v\|\left(\|e_{k_\delta-1}\|+\|e_{k-1}\|\right).
\end{eqnarray}
In order to estimate $J$, we write
\begin{equation*}
J=J_1+J_2+J_3,
\end{equation*}
where 
\begin{eqnarray*}
J_1&:=\left(1-\frac{\alpha_{k-1}}{\alpha_{k_\delta-1}}\right)(\alpha_{k-1}I+\A)^{-1} F'(x^\dag)^*(F(x_{k_\delta})-y),\\
J_2&:=\left(1-\frac{\alpha_{k-1}}{\alpha_{k_\delta-1}}\right) (\alpha_{k-1} I+\A)^{-1} F'(x^\dag)^*\left[F'(x^\dag) e_{k_\delta}-F(x_{k_\delta})+y\right],\\
J_3&:=\left(1-\frac{\alpha_{k-1}}{\alpha_{k_\delta-1}}\right)(\alpha_{k-1} I+\A)^{-1} \A\left[\alpha_{k_\delta-1}(\alpha_{k_\delta-1} I+\A)^{-1} e_0-e_{k_\delta}\right].
\end{eqnarray*}
By using the argument in the proof of \cite[Lemma 4.4]{Jin01} we can see that
\begin{equation*}
\|J_1\|\lesssim \frac{1}{\sqrt{\alpha_k}}\|\alpha_{k_\delta}^{1/2}(\alpha_{k_\delta} I+\B)^{-1/2}(F(x_{k_\delta})-y)\|.
\end{equation*}
Also, by using (\ref{1}) and noting $\alpha_{k-1}\le \alpha_{k_\delta-1}$, it is easy to see that
\begin{equation*}
\|J_2\|\lesssim \frac{1}{\sqrt{\alpha_{k-1}}} L\|e_{k_\delta}\|^2.
\end{equation*}
From Proposition \ref{P1} and Lemma \ref{L5} we have
\begin{equation*}
\|e_{k_\delta}\|\le \|e_{k_\delta}^\delta\|+\|x_{k_\delta}^\delta-x_{k_\delta}\|\lesssim \|v\|^{1/2} \delta^{1/2}+\frac{\delta}{\sqrt{\alpha_{k_\delta}}}.
\end{equation*}
Recall that $k_\delta\le \tilde{k}_\delta$ which implies $\frac{\delta}{\sqrt{\alpha_{k_\delta}}}
\lesssim \|v\|^{1/2}\delta^{1/2}$. Thus
\begin{equation*}
\|e_{k_\delta}\|\lesssim \|v\|^{1/2}\delta^{1/2}.
\end{equation*}
Since $k\le \tilde{k}_\delta$, we then obtain
\begin{equation*}
\frac{\|e_{k_\delta}\|}{\sqrt{\alpha_{k-1}}}
\lesssim\frac{1}{\sqrt{\alpha_{k-1}}}\|v\|^{1/2}\delta^{1/2}\lesssim \|v\|.
\end{equation*}
Therefore 
\begin{equation*}
\|J_2\|\lesssim L\|v\|\|e_{k_\delta}\|. 
\end{equation*}
By using (\ref{70}), $J_3$ can be estimated as 
\begin{equation*}
\|J_3\|\le \|\alpha_{k_\delta-1}(\alpha_{k_\delta-1} I+\A)^{-1} e_0-e_{k_\delta}\|\lesssim L\|v\|\|e_{k_\delta-1}\|.
\end{equation*}
Combining the above estimates on $J_1$, $J_2$ and $J_3$, we obtain
\begin{equation*}
\|J\|\lesssim \frac{\|\alpha_{k_\delta}^{1/2}(\alpha_{k_\delta} I+\B)^{-1/2}(F(x_{k_\delta})-y)\|}
{\sqrt{\alpha_k}}
+L\|v\|\left(\|e_{k_\delta}\|+\|e_{k_\delta-1}\|\right).
\end{equation*}
This together with (\ref{3.10}) and (\ref{3.100}) gives
\begin{equation*}
\fl \|x_{k_\delta}-x_k\|\lesssim \frac{\|\alpha_{k_\delta}^{1/2}
(\alpha_{k_\delta} I+\B)^{-1/2}(F(x_{k_\delta})-y)\|}{\sqrt{\alpha_k}}
+L\|v\|\left(\|e_k\|+\|e_{k_\delta}\|\right).
\end{equation*}
By the smallness of $L\|v\|$, we thus conclude the proof.
\hfill $\Box$

\section{\bf Proof of Theorem \ref{T1}}

In this section we will complete the proof of the main result, Theorem \ref{T1}. 
In order to apply Proposition \ref{P2}, we need the following estimates.

\begin{lemma}\label{L7}
Assume that (\ref{81}), (\ref{1}) and (\ref{94}) hold and that $x_0-x^\dag=F'(x^\dag)^* v$ for some $v\in 
{\mathcal N}(F'(x^\dag)^*)^\perp$. Let $k_\delta$ be the integer determined by (\ref{4}) 
with $\tau>2$. If $L\|v\|$ is sufficiently small, then we have
\begin{equation*}
\|\alpha_{k_\delta}(\alpha_{k_\delta} I+\B)^{-1/2} (F(x_{k_\delta})-y)\|\lesssim \delta
\end{equation*}
and
\begin{equation*}
\delta\lesssim \|\alpha_k^{1/2}(\alpha_k I+\B)^{-1/2} (F(x_k)-y)\|
\end{equation*}
for all $0\le k<k_\delta$. 
\end{lemma}

\noindent{\bf Proof.}
For $0\le k\le k_\delta$ we set
\begin{eqnarray*}
a_k&:=\|\alpha_k^{1/2}(\alpha_k I+\B)^{-1/2} (F(x_k)-y)\|^2,\\
b_k&:=\|\alpha_k^{1/2}(\alpha_k I+\B_k^\delta)^{-1/2}(F(x_k)-y)\|^2.
\end{eqnarray*}
It then follows from \cite[Proposition 3.4]{TJ} and Lemma \ref{L3} that
\begin{equation*}
|a_k-b_k|\lesssim \frac{1}{\sqrt{\alpha_k}} L\|e_k^\delta\|(a_k+b_k)\lesssim L\|v\|(a_k+b_k).
\end{equation*}
By the smallness of $L\|v\|$, we have $a_k\approx b_k$. Thus it suffices to show that
\begin{equation*}
\sqrt{b_{k_\delta}}\lesssim \delta \quad \mbox{and} \quad \delta\lesssim \sqrt{b_k} \quad \mbox{for } 0\le k<k_\delta.
\end{equation*}
By using (\ref{1}), Lemma \ref{L5} and (\ref{3.4}) we have for $0\le k<k_\delta$
\begin{eqnarray*}
\fl \sqrt{b_k}\ge \|\alpha_k^{1/2}(\alpha_k I+\B_k^\delta)^{-1/2}(F(x_k^\delta)-y^\delta)\|
-\|F(x_k)-F(x_k^\delta)-F'(x_k^\delta)(x_k-x_k^\delta)\|\\
 -\|\alpha_k^{1/2} (\alpha_k I+\B_k^\delta)^{-1/2} F'(x_k^\delta)(x_k-x_k^\delta)\|-\delta\\
\ge (\tau-1)\delta -\alpha_k^{1/2} \|x_k-x_k^\delta\|-\frac{1}{2} L\|x_k-x_k^\delta\|^2\\
\ge (\tau-2)\delta- CL\frac{\delta^2}{\alpha_k}\\
\ge (\tau-2-CL\|v\|)\delta,
\end{eqnarray*}
where $C$ is a universal constant independent of $\delta$. 
Using $\tau>2$ and the smallness of $L\|v\|$ we obtain
$\delta\lesssim \sqrt{b_k}$ for all $0\le k<k_\delta$.

Similarly, we have
\begin{eqnarray*}
\sqrt{b_{k_\delta}}&\le (\tau+1)\delta +\alpha_{k_\delta}^{1/2}\|x_{k_\delta}-x_{k_\delta}^\delta\|
+\frac{1}{2}L\|x_{k_\delta}^\delta-x_{k_\delta}\|^2\\
&\lesssim \delta+L\frac{\delta^2}{\alpha_{k_\delta}}\\
&\lesssim \delta.
\end{eqnarray*}
The proof is thus complete.
\hfill $\Box$\\

\noindent{\bf Proof of Theorem \ref{T1}.}
We first prove (\ref{0.8}). 
Note that for $k>\tilde{k}_\delta$, we have $\frac{\delta}{\sqrt{\alpha_k}}\ge c_0\alpha_k^{1/2} \|v\|$, 
while Lemma \ref{L3} implies that $\|e_{k_\delta}^\delta\|\lesssim \alpha_k^{1/2}\|v\|$. 
Therefore, in order to complete the proof, it suffices to show
\begin{equation*}
\|e_{k_\delta}^\delta\|\lesssim \inf\left\{\|e_k\|+\frac{\delta}{\sqrt{\alpha_k}}: 0\le k\le \tilde{k}_\delta\right\}.
\end{equation*}
Note that for $k_\delta\le k\le \tilde{k}_\delta$, we have from Lemma \ref{L5}, Proposition \ref{P2} and Lemma \ref{L7} that
\begin{eqnarray*}
\|e_{k_\delta}^\delta\|&\lesssim \|e_{k_\delta}\|+\frac{\delta}{\sqrt{\alpha_{k_\delta}}}\\
&\lesssim \|e_k\|+\frac{\|\alpha_{k_\delta}^{1/2}(\alpha_{k_\delta} I+\B)^{-1/2} (F(x_{k_\delta})-y)\|}{\sqrt{\alpha_k}}
+\frac{\delta}{\sqrt{\alpha_k}}\\
&\lesssim \|e_k\|+\frac{\delta}{\sqrt{\alpha_k}},
\end{eqnarray*}
while for $0\le k<k_\delta$ we have from Lemma \ref{L2} and Lemma \ref{L7} that
\begin{eqnarray*}
\fl \|e_{k_\delta}^\delta\|\lesssim \|e_{k_\delta}\|+\frac{\delta}{\sqrt{\alpha_{k_\delta}}}\\
\fl \qquad\quad \lesssim \|e_k\|+\frac{1}{\sqrt{\alpha_{k_\delta}}}\|\alpha_{k_\delta-1}^{1/2}(\alpha_{k_\delta-1} I+\B)^{-1/2} (F(x_{k_\delta-1})-y)\|\\
\fl\qquad\quad\lesssim  \|e_k\|+\|(\alpha_{k_\delta-1} I+\B)^{-1/2} (F(x_{k_\delta-1})-y)\|\\
\fl \qquad\quad\lesssim  \|e_k\|+\|(\alpha_{k_\delta-1} I+\B)^{-1/2} F'(x^\dag) e_{k_\delta-1}\|
+\frac{\|F(x_{k_\delta-1})-y-F'(x^\dag) e_{k_\delta-1}\|}{\sqrt{\alpha_{k_\delta-1}}}\\
\fl \qquad\quad\lesssim  \|e_k\|+\|e_{k_\delta-1}\|+\frac{1}{\sqrt{\alpha_{k_\delta-1}}} L\|e_{k_\delta-1}\|^2\\
\fl \qquad\quad\lesssim \|e_k\|+L\|v\|\|e_{k_\delta-1}\|\\
\fl \qquad\quad \lesssim \|e_k\|.
\end{eqnarray*}
The proof of (\ref{0.8}) is complete.

Next we prove (\ref{0.8.5}). Note that $\beta_k\le \alpha_k^{\nu/2}\|\omega\|$
under the condition on $x_0-x^\dag$. We have from Lemma \ref{L2} that $\|e_k\|\lesssim
\alpha_k^{\nu/2}\|\omega\|$. Thus it follows from (\ref{0.8}) that 
\begin{equation*}
\|e_{k_\delta}^\delta\|\lesssim \inf\left\{\alpha_k^{\nu/2}\|\omega\|+\frac{\delta}{\sqrt{\alpha_k}}: k=0, 1, \cdots\right\}.
\end{equation*}
Now we introduce the integer $\bar{k}_\delta$ such that
\begin{equation*}
\alpha_{\bar{k}_\delta}\le \left(\frac{\delta}{\|\omega\|}\right)^{2/(1+\nu)}<\alpha_k, \quad 0\le k<\bar{k}_\delta.
\end{equation*}
Then it is readily to see that
\begin{equation*}
\|e_{k_\delta}^\delta\|\lesssim \alpha_{\bar{k}_\delta}^{\nu/2}\|\omega\|+\frac{\delta}{\sqrt{\alpha_{\bar{k}_\delta}}}
\lesssim \|\omega\|^{1/(1+\nu)}\delta^{\nu/(1+\nu)}.
\end{equation*}
\hfill $\Box$

\end{document}